\documentclass{amsart}

\usepackage{fancyhdr}
\usepackage{lastpage}
\usepackage{stmaryrd,yhmath}

\newcommand{\bdis}{\begin{displaymath}}
\newcommand{\edis}{\end{displaymath}}
\newcommand{\be}{\begin{equation}}
\newcommand{\ee}{\end{equation}}
\newcommand{\mbb}{\mathbb}
\newcommand{\mcal}{\mathcal}

\newcommand{\vp}{\varphi}

\newcommand{\mT}{\mathring{T}}

\newcommand{\zf}{\zeta\left(\frac{1}{2}+it\right)}

\theoremstyle{definition}

\newtheorem{cor}[]{Corollary}

\theoremstyle{remark}
\newtheorem{remark}[]{Remark}

\newtheorem*{mydef1}{{\bf Theorem}}

\numberwithin{equation}{section}

\begin{document}

\title{Jacob's ladders and new orthogonal systems generated by Jacobi polynomials}

\author{Jan Moser}

\address{Department of Mathematical Analysis and Numerical Mathematics, Comenius University, Mlynska Dolina M105, 842 48 Bratislava, SLOVAKIA}

\email{jan.mozer@fmph.uniba.sk}


\begin{abstract}
Is is shown in this paper that there is a connection between the Riemann zeta-function $\zf$ and the classical Jacobi's polynomials, i.e. the Legendre 
polynomials, Chebyshev polynomials of the first and the second kind, \dots . 
\end{abstract}

\maketitle

\section{The result} 

\subsection{}  

In this paper we obtain some new properties of the signal 
\bdis 
Z(t)=e^{i\vartheta(t)}\zf 
\edis  
that is generated by the Riemann zeta-function, where 
\bdis 
\vartheta(t)=-\frac t2\ln\pi+\text{Im}\ln\Gamma\left(\frac 14+i\frac t2\right)=\frac t2\ln\frac{t}{2\pi}-\frac t2-\frac \pi 8+\mcal{O}\left(\frac 1t\right). 
\edis  
Let us remind that 
\bdis 
\tilde{Z}^2(t)=\frac{{\rm d}\vp_1(t)}{{\rm d}t},\ \vp_1(t)=\frac 12\vp(t) 
\edis  
where 
\be \label{1.1} 
\tilde{Z}^2(t)=\frac{Z^2(t)}{2\Phi^\prime_\vp[\vp(t)]}=\frac{\left|\zf\right|^2}{\left\{ 1+\mcal{O}\left(\frac{\ln\ln t}{\ln t}\right)\right\}\ln t}
\ee  
(see \cite{1}, (3.9); \cite{3}, (1.3); \cite{7}, (1.1), (3.1), (3.2)), and $\vp(t)$ is the Jacob's ladder, i.e. a solution of the nonlinear integral 
equation (see \cite{1}) 
\bdis 
\int_0^{\mu[X(T)]}Z^2(t)e^{-\frac{2}{X(T)}t}{\rm d}t=\int_0^T Z^2(t){\rm d}t . 
\edis  

\subsection{} 

The system of the Jacobi polynomials 
\bdis 
\begin{split} 
& (1-x)^\alpha(1+x)^\beta P^{(\alpha,\beta)}_n(x)=\frac{(-1)^n}{2^nn!}\frac{{\rm d}^n}{{\rm d}x^n}\left[(1-x)^{\alpha+n}(1+x)^{\beta+n}\right], \\ 
& x\in [-1,1],\ n=0,1,2,\dots ,\ \alpha,\beta>-1 
\end{split} 
\edis  
is the well-known system of orthogonal polynomials on the segment $x\in [-1,1]$ with the weight function 
\bdis 
(1-x)^\alpha(1+x)^\beta , 
\edis   
i.e. the following formulae hold true (comp. \cite{18}) 
\be\label{1.2} 
\begin{split} 
& \int_{-1}^1 (1-x)^\alpha(1+x)^\beta P^{(\alpha,\beta)}_m(x)P^{(\alpha,\beta)}_n(x){\rm d}x=0,\ m\not=n , \\ 
& \int_{-1}^1 (1-x)^\alpha(1+x)^\beta \left[P^{(\alpha,\beta)}_n(x)\right]^2{\rm d}x=
\frac{2^{\alpha+\beta+1}}{2n+\alpha+\beta+1}\frac{\Gamma(n+\alpha+1)\Gamma(n+\beta+1)}{n!\Gamma(n+\alpha+\beta+1)} . 
\end{split}
\ee 

It is shown in this paper that the $\tilde{Z}^2$-transformation of the Jacobi's polynomials generates a new system of orthogonal functions connected 
with $|\zf|^2$. In this direction, the following theorem holds true. 

\begin{mydef1} 
 Let $x=t-T-1,\ t\in [T,T+2]$ and 
\bdis 
\vp_1\{[\mT,\widering{T+2}]\}=[T,T+2],\ T\geq T_0[\vp_1]. 
\edis  
Then the system of functions 
\bdis 
P^{(\alpha,\beta)}_n(\vp_1(t)-T-1),\ t\in[\mT,\widering{T+2}],\ n=0,1,2,\dots 
\edis 
is the orthogonal system on $[\mT,\widering{T+2}]$ with the weight function given by 
\bdis 
(T+2-\vp_1(t))^\alpha(\vp_1(t)-T)^\beta\tilde{Z}^2(t) , 
\edis  
i.e. the following system of new-type integrals 
\be \label{1.3} 
\begin{split} 
& \int_{\mT} ^{\widering{T+2}}P^{(\alpha,\beta)}_m(\vp_1(t)-T-1)P^{(\alpha,\beta)}_n(\vp_1(t)-T-1) \\ 
& (T+2-\vp_1(t))^\alpha(\vp_1(t)-T)^\beta\tilde{Z}^2(t){\rm d}t=0, \ m\not=n , \\ 
& \int_{\mT} ^{\widering{T+2}}\left[ P^{(\alpha,\beta)}_n(\vp_1(t)-T-1)\right]^2(T+2-\vp_1(t))^\alpha(\vp_1(t)-T)^\beta\tilde{Z}^2(t){\rm d}t= \\ 
& =\frac{2^{\alpha+\beta+1}}{2n+\alpha+\beta+1}\frac{\Gamma(n+\alpha+1)\Gamma(n+\beta+1)}{n!\Gamma(n+\alpha+\beta+1)} , \\ 
& m,n=0,1,2,\dots , \ \text{for all} \ \mT\geq \vp_1^{-1}(T),\ T\geq T_0[\vp_1]
\end{split}
\ee  
is obtained, where $\vp_1(t)-T-1\in [-1,1]$, and 
\be \label{1.4} 
\rho\{[-1,1];[\mT,\widering{T+2}]\}\sim T, \ T\to\infty . 
\ee 
\end{mydef1}

\begin{remark} 
This Theorem gives the contact point between the functions $\zf$, $\vp_1(t)$ and the Jacobi polynomials $P^{(\alpha,\beta)}_n(x)$. 
\end{remark} 

\subsection{} 

The seconf formula in (\ref{1.3}) via the mean-value theorem (comp. (\ref{1.1}) leads to 

\begin{cor} 
\be \label{1.5} 
\begin{split} 
& \int_{\mT} ^{\widering{T+2}}\left[ P^{(\alpha,\beta)}_n(\vp_1(t)-T-1)\right]^2(T+2-\vp_1(t))^\alpha(\vp_1(t)-T)^\beta
\left|\zf\right|^2{\rm d}t\sim \\ 
& \sim \frac{2^{\alpha+\beta+1}}{2n+\alpha+\beta+1}\frac{\Gamma(n+\alpha+1)\Gamma(n+\beta+1)}{n!\Gamma(n+\alpha+\beta+1)}\ln \mT,\ \mT\to\infty \\ 
& n=0,1,2,\dots . 
\end{split}
\ee 
\end{cor} 

This paper is a continuation of the series \cite{1}-\cite{17}. 

\section{Orthogonal systems generated by Legendre polynomials} 

If $\alpha=\beta=0$ then $P^{(\alpha,\beta)}_m(x)=P_n(x)$ is the Legendre polynomial. In this case, our Theorem implies the following 

\begin{cor} 
The system of functions 
\bdis 
P_n(\vp_1(t)-T-1),\ t\in [\mT,\widering{T+2}],\ n=0,1,2,\dots 
\edis  
is the orthogonal system on the segment $[\mT,\widering{T+2}]$ with the weight function $\tilde{Z}^2(t)$, i.e. the following system of new-type integrals 
\be \label{2.1} 
\begin{split}
& \int_{\mT} ^{\widering{T+2}} P_m(\vp_1(t)-T-1)P_n(\vp_1(t)-T-1)\tilde{Z}^2(t){\rm d}t=0,\ m\not=n, \\ 
& \int_{\mT} ^{\widering{T+2}} \left[ P_n(\vp_1(t)-T-1)\right]^2\tilde{Z}^2(t){\rm d}t=\frac{2}{2n+1} , \\ 
& m,n=0,1,2,\dots ,\ \text{for all} \ \mT\geq \vp_1^{-1}(T),\ T\geq T_0[\vp] 
\end{split} 
\ee  
is obtained. 
\end{cor}

From the second formula in (\ref{2.1}) we obtain (comp. (\ref{1.5})) 

\begin{cor} 
\be\label{2.2} 
\begin{split} 
& \int_{\mT} ^{\widering{T+2}}\left[ P_n(\vp_1(t)-T-1)\right]^2\left|\zf\right|^2{\rm d}t\sim \frac{2}{2n+1}\ln\mT,\ \mT\to\infty,\\ 
& n=0,1,2,\dots . 
\end{split} 
\ee 
\end{cor}

\section{Orthogonal systems generated by Chebyshev polynomials of the first and the second kind} 

\subsection{} 

If $\alpha=\beta=-\frac 12$ then 
\bdis 
P^{(-\frac 12,-\frac 12)}_n(x)=\frac{(2n-1)!!}{(2n)!!}T_n(x),\ n=0,1,2,\dots 
\edis  
where $T_n(x)$ is the Chebyshev polynomial of the first kind. In this case we obtain from our Theorem 

\begin{cor} 
The system of functions 
\bdis 
T_n(\vp_1(t)-T-1),\ t\in [\mT,\widering{T+2}],\ n=0,1,2,\dots 
\edis  
is the orthogonal system of functions with the weight function 
\bdis 
\frac{\tilde{Z}^2(t)}{\sqrt{1-(\vp_1(t)-T-1)^2}} , 
\edis   
i.e. the following system of the new-type integrals 
\be \label{3.1} 
\begin{split} 
& \int_{\mT}^{\widering{T+2}}T_m(\vp_1(t)-T-1)T_n(\vp_1(t)-T-1)\frac{\tilde{Z}^2(t)}{\sqrt{1-(\vp_1(t)-T-1)^2}}{\rm d}t=0,\ m\not= n, \\ 
& \int_{\mT}^{\widering{T+2}}\left[ T_n(\vp_1(t)-T-1)\right]^2\frac{\tilde{Z}^2(t)}{\sqrt{1-(\vp_1(t)-T-1)^2}}{\rm d}t=
\left\{\begin{array}{rcl} \frac{\pi}{2} & , & n\geq 1, \\ \pi & , & n=0, \end{array} \right. 
\end{split}
\ee 
$m,n=0,1,2,\dots $, for all $\mT\geq \vp_1^{-1}(T),\ T\geq T_0[\vp_1]$, is obtained. 
\end{cor}

From the second formula in (\ref{3.1}) we obtain (comp. (\ref{1.5})) 

\begin{cor} 
\be \label{3.2} 
\begin{split} 
& \int_{\mT}^{\widering{T+2}}\left[ T_n(\vp_1(t)-T-1)\right]^2\frac{\left|\zf\right|^2}{\sqrt{1-(\vp_1(t)-T-1)^2}}{\rm d}t\sim \\ 
& \sim \left\{\begin{array}{rcl} \frac{\pi}{2}\ln\mT & , & n\geq 1 , \\ \pi\ln\mT & , & n=0, \end{array} \right. \ \mT\to\infty , \ n=0,1,2,\dots \ . 
\end{split} 
\ee 
\end{cor}

\begin{remark} 
From (\ref{3.2}) (see the second formula; since $T_0(x)=1$) the canonical asymptotic formula 
\be \label{3.3} 
\int_{\mT}^{\widering{T+2}}\frac{\left|\zf\right|^2}{\sqrt{1-(\vp_1(t)-T-1)^2}}{\rm d}t\sim \pi\ln\mT,\ \mT\to\infty 
\ee  
follows. 
\end{remark} 

\subsection{} 

If $\alpha=\beta=\frac 12$ then 
\bdis 
P^{(\frac 12,\frac 12)}_n(x)=2\frac{(2n-1)!!}{(2n+2)!!}U_n(x),\ n=0,1,2,\dots 
\edis  
where $U_n(x)$ is the Chebyshev polynomial of the second kind. Then, from our Theorem, we obtain 
\begin{cor} 
The system of the functions 
\bdis 
U_n(\vp_1(t)-T-1),\ t\in [\mT,\widering{T+2}],\ n=0,1,2,\dots 
\edis  
is the orthogonal system with the weight function 
\bdis 
\sqrt{1-(\vp_1(t)-T-1)^2}\tilde{Z}^2(t) , 
\edis  
i.e. the following system of the new-type integrals 
\be \label{3.4} 
\begin{split} 
& \int_{\mT}^{\widering{T+2}}U_m(\vp_1(t)-T-1)U_n(\vp_1(t)-T-1)\\ 
& \sqrt{1-(\vp_1(t)-T-1)^2}\tilde{Z}^2(t){\rm d}t=0,\ m\not= n, \\ 
& \int_{\mT}^{\widering{T+2}}\left[ U_n(\vp_1(t)-T-1)\right]^2\sqrt{1-(\vp_1(t)-T-1)^2}\tilde{Z}^2(t){\rm d}t=\frac{\pi}{2}, 
\end{split} 
\ee 
$m,n=0,1,2,\dots $, for all $\mT\geq\vp_1^{-1}(T),\ T\geq T_0[\vp_1]$ is obtained. 
\end{cor}

From the second formula in (\ref{3.4}) we obtain (comp. (\ref{1.5})) 

\begin{cor} 
\be \label{3.5} 
\begin{split} 
& \int_{\mT}^{\widering{T+2}}\left[ U_n(\vp_1(t)-T-1)\right]^2\sqrt{1-(\vp_1(t)-T-1)^2}\left|\zf\right|^2{\rm d}t\sim \\ 
& \sim \frac{\pi}{2}\ln\mT,\ n=0,1,2,\dots ,\ \mT\to\infty . 
\end{split}
\ee 
\end{cor} 

\begin{remark} 
From (\ref{3.5}) (see the second formula; since $U_0(x)=1$) the canonical asymptotic formula 
\be \label{3.6} 
\int_{\mT}^{\widering{T+2}}\sqrt{1-(\vp_1(t)-T-1)^2}\left|\zf\right|^2{\rm d}t\sim \frac{\pi}{2}\ln\mT,\ \mT\to\infty 
\ee  
follows. 
\end{remark} 

\begin{remark} 
The Riemann zeta-function $\zf$ is connected with the classical orthogonal polynomials of Legendre and of Chebyshev by formulae (\ref{2.1}), (\ref{2.2}), 
(\ref{3.1})-(\ref{3.6}), respectivelly. 
\end{remark} 

\section{Proof of Theorem} 

\subsection{} 

Let us remind that the following lemma holds true (see \cite{6}, (2.5); \cite{7}, (3.3)): for every integrable function (in the Lebeague sense) 
$f(x), x\in[\vp_1(T),\vp_1(T+U)]$ we have 
\be \label{4.1} 
\int_T^{T+U}f[\vp_1(t)]\tilde{Z}^2(t){\rm d}t=\int_{\vp_1(T)}^{\vp_1(T+U)}f(x){\rm d}x,\ U\in (0,T/\ln T] 
\ee  
where 
\be \label{4.2} 
t-\vp_1(t)\sim (1-c)\pi(t) , 
\ee  
$c$ is the Euler's constant and $\pi(t)$ is the prime-counting function. In the case (comp. Theorem) $T=\vp_1(\mT), T+U=\vp_1(\widering{T+U})$, we 
obtain from (\ref{4.1}) 
\be \label{4.3} 
\int_{\mT}^{\widering{T+U}}f[\vp_1(t)]\tilde{Z}^2(t){\rm d}t=\int_T^{T+U}f(x){\rm d}x . 
\ee  

\subsection{} 

Putting 
\bdis 
\begin{split} 
& f(t)=P^{(\alpha,\beta)}_m(t-T-1)P^{(\alpha,\beta)}_n(t-T-1)(T+2-t)^\alpha(t-T)^\beta;\ U=2, 
\end{split}  
\edis  
we have by (\ref{4.3}) and (\ref{1.2}) the following $\tilde{Z}^2$-transformation 
\bdis 
\begin{split} 
& \int_{\mT}^{\widering{T+U}}P^{(\alpha,\beta)}_m(\vp_1(t)-T-1)P^{(\alpha,\beta)}_n(\vp_1(t)-T-1) \\ 
& (T+2-\vp_1(t))^\alpha(\vp_1(t)-T)^\beta\tilde{Z}^2(t){\rm d}t= \\ 
& \int_T^{T+2} P^{(\alpha,\beta)}_m(t-T-1)P^{(\alpha,\beta)}_n(t-T-1)(T+2-t)^\alpha(t-T)^\beta{\rm d}t= \\ 
& \int_{-1}^1 P^{(\alpha,\beta)}_m(x)P^{(\alpha,\beta)}_n(x)(1-x)^\alpha(1+x)^\beta{\rm d}x=0,\ m\not= n, 
\end{split} 
\edis  
where $t=x+T+1$, i.e. the first formula in (\ref{1.3}) holds true. Similarly we obtain the second formula in (\ref{1.3}). Since (\ref{4.2}) implies 
$\mT\to T$ then (\ref{1.4}) holds true. 

\subsection{} 

Next, from (\ref{4.2}) we obtain 
\bdis 
\begin{split} 
& \mT-\vp_1(\mT)=\mT-T=\mcal{O}\left(\frac{\mT}{\ln \mT}\right) , \\ 
& \widering{T+2}-\vp_1(\widering{T+2})=\widering{T+2}-T-2=\mcal{O}\left(\frac{\mT}{\ln \mT}\right), 
\end{split} 
\edis 
and subsequently 
\bdis 
\widering{T+2}-\mT=\mcal{O}\left(\frac{\mT}{\ln \mT}\right) , 
\edis  
and for $\xi\in(\mT,\widering{T+2})$ we have 
\be \label{4.4} 
\ln\xi=\ln\mT+\mcal{O}\left(\frac{\widering{T+2}-\mT}{\mT}\right)=\ln \mT+\mcal{O}\left(\frac{1}{\ln\mT}\right) . 
\ee  
The property (\ref{4.4}) was used in (\ref{1.5}),(\ref{2.2}),(\ref{3.2}),(\ref{3.5}).

\thanks{I would like to thank Michal Demetrian for helping me with the electronic version of this work.}


\begin{thebibliography}{29}% Replace 9 by 99 if 10 or more references
%
% Please note the use of "\and" between author names below
\bibitem{1}
J. Moser, `Jacob's ladders and the almost exact asymptotic representation of the Hardy-Littlewood integral', Math. Notes 2010, {\bf 88}, 
pp. 414-422, arXiv:0901.3973.
%
\bibitem{2}
J. Moser, `Jacob's ladders and the tangent law for short parts of the Hardy-Littlewood integral', (2009), arXiv:0906.0659.
%
\bibitem{3}
J. Moser, `Jacob's ladders and the multiplicative asymptotic formula for short and microscopic parts of the Hardy-Littlewood integral', (2009),
arXiv:0907.0301.
%
\bibitem{4}
J. Moser, `Jacob's ladders and the quantization of the Hardy-Littlewood integral', (2009), arXiv:0909.3928.
%
\bibitem{5}
J. Moser,
`Jacob's ladders and the first asymptotic formula for the expression of the sixth order $|\zeta(1/2+i\varphi(t)/2)|^4|\zeta(1/2+it)|^2$', (2009),
arXiv:0911.1246.
%
\bibitem{6}
J. Moser,
`Jacob's ladders and the first asymptotic formula for the expression of the fifth order
$Z[\varphi(t)/2+\rho_1]Z[\varphi(t)/2+\rho_2]Z[\varphi(t)/2+\rho_3]\hat{Z}^2(t)$ for the collection of disconnected sets`, (2009),
arXiv:0912.0130.
%
\bibitem{7}
J. Moser, `Jacob's ladders, the iterations of Jacob's ladder $\vp_1^k(t)$ and asymptotic formulae for the integrals of the products
$Z^2[\varphi^n_1(t)]Z^2[\varphi^{n-1}(t)]\cdots Z^2[\varphi^0_1(t)]$ for arbitrary fixed $n\in \mbb{N}$` (2010), arXiv:1001.1632.
%
\bibitem{8}
J. Moser, `Jacob's ladders and the asymptotic formula for the integral of the eight order expression
$|\zeta(1/2+i\vp_2(t))|^4|\zeta(1/2+it)|^4$`, (2010), arXiv:1001.2114.
%
\bibitem{9}
J. Moser, `Jacob's ladders and the asymptotically approximate solutions of a nonlinear diophantine equation`, (2010), arXiv: 1001.3019.
%
\bibitem{10}
J. Moser, `Jacob's ladders and the asymptotic formula for short and microscopic parts of the Hardy-Littlewood integral of the function
$|\zeta(1/2+it)|^4$`, (2010), arXiv:1001.4007.
%
\bibitem{11}
J. Moser, `Jacob's ladders and the nonlocal interaction of the function $|\zeta(1/2+it)|$ with $\arg\zeta(1/2+it)$ on the distance
$\sim (1-c)\pi(t)$`, (2010), arXiv: 1004.0169.
%
\bibitem{12}
J. Moser, `Jacob's ladders and the $\tilde{Z}^2$ - transformation of polynomials in $\ln\vp_1(t)$`, (2010), arXiv: 1005.2052.
%
\bibitem{13}
J. Moser,
`Jacob's ladders and the oscillations of the function $|\zf|^2$ around the main part of its mean-value; law of the almost exact equality of the
corresponding areas`, (2010), arXiv: 1006.4316
%
\bibitem{14}
J. Moser,
`Jacob's ladders and the nonlocal interaction of the function $Z(t)$ with the function $\tilde{Z}^2(t)$ on the distance $\sim (1-c)\pi(t)$ for a
collection of disconneted sets`, (2010), arXiv: 1006.5158
%
\bibitem{15}
J. Moser, `Jacob's ladders and the $\tilde{Z}^2$-transformation of the orthogonal system of trigonometric functions`, (2010), arXiv: 1007.0108.
%
\bibitem{16}
J. Moser, `Jacob's ladders and the nonlocal interaction of the function $Z^2(t)$ with the function $\tilde{Z}^2(t)$ on the distance $\sim (1-c)\pi(t)$
for the collections of disconnected sets`, (2010), arXiv: 1007.5147.
%
\bibitem{17}
J. Moser, `Jacob's ladders and some new consequences from A. Sleberg's formula`, (2010), arXiv: 1010.0868.
%
\bibitem{18}
G. Szeg\" o, `Orthogonal polynomials`, New York, 1959. 

\end{thebibliography}
\end{document}